\documentclass{article}
\usepackage{amsmath,epsfig,latexsym,amssymb,cite,graphicx,a4}
\newlength{\indentedwidth} \newdimen\mathindent
\indentedwidth=\mathindent

\newtheorem{proposition}{Proposition}[section]

\newtheorem{definition}{Definition}[section]

\usepackage{epsf,graphics}
\DeclareMathAlphabet{\mathpzc}{OT1}{pzc}{m}{it}
\begin{document}
\vskip 0.5cm
\begin{center}
{\Large \bf Lax Operator for the Quantised Orthosymplectic
Superalgebra $U_q[osp(2|n)]$}
\end{center}
\vskip 0.8cm
\centerline{K.A. Dancer,%
\footnote{\tt dancer@maths.uq.edu.au} M.D. Gould%
\footnote{\tt mdg@maths.uq.edu.au} and J. Links%
\footnote{\tt jrl@maths.uq.edu.au}}
\vskip 0.9cm \centerline{\sl\small Centre for Mathematical Physics,
School of Physical Sciences, } \centerline{\sl\small The University of
Queensland, Brisbane 4072, Australia.}
\vskip 0.9cm
\begin{abstract} Each quantum superalgebra is a quasi-triangular Hopf
  superalgebra, so contains a \textit{universal $R$-matrix} in the tensor 
  product algebra which satisfies
  the Yang-Baxter equation.  Applying the vector representation $\pi$,
  which acts on the vector module $V$, to one side of a universal
  $R$-matrix gives a Lax operator.  In this paper a Lax operator is
  constructed for the $C$-type quantum superalgebras $U_q[osp(2|n)]$. 
  This can in turn be used to find a solution to
  the Yang-Baxter equation acting on $V \otimes V \otimes W$ where $W$
  is an arbitrary $U_q[osp(2|n)]$ module. The case $W=V$ is included here
  as an example.  \end{abstract}


\section{Introduction}

The representation theory of quantum superalgebras is a 
useful tool for modelling supersymmetric 
systems, especially in cases where the system is integrable.  
Solutions of the spectral parameter dependent 
Yang--Baxter equation provide the foundation for constructing integrable models,
and also facilitate the associated Bethe ansatz method for deriving the exact 
solution. Examples of integrable systems obtained in this manner include
those describing strongly correlated electrons on a one-dimensional
lattice \cite{Foerster,Gonzalez,uq,uqagain,Ramos} and supersymmetric field theories 
\cite{Saleur1,Saleur2,Saleur3}. 

A means to construct algebraic solutions to the Yang--Baxter equation is given 
by the supersymmetric version of the
Quantum Double construction to obtain a universal $R$-matrix $\mathcal{R}$ 
\cite{gzb}. For the case of an affine quantum superalgebra,
loop representations of $\mathcal{R}$ then yield solutions of the Yang--Baxter
equation with spectral parameter. An alternate method is to start with a 
representation of the universal 
$R$-matrix for the 
non-affine case, and then introduce a spectral parameter by using tensor 
product graph methods \cite{dglz, GouldZhang00}.  
Even in the non-affine case, however, determining the
universal $R$-matrix is a difficult challenge. The known explicit matrix
solutions have been constructed by symmetry considerations rather than 
resulting directly from representations of the universal $R$-matrix.
One way to simplify the problem, whilst maintaininng a level of generality,
is to instead construct the Lax operator 
$R = (\pi \otimes I) \mathcal{R}$ where $\pi$ denotes the vector 
representation of the quantum superalgebra.  This has been done for the $A$, 
$B$ and $D$ families of quantum superalgebras in \cite{Zhang2,us}, but the 
solution does not automatically extend to the $C$-series superalgebras 
$U_q[osp(2|n)]$ because of differences in the root systems.  Here the same 
methods used in \cite{us} are adapted to construct the Lax 
operator for $U_q[osp(2|n)]$, thus completing the calculation of Lax operators 
for all non-exceptional quantum superalgebras.   

We begin in Section 2 with the construction of $U_q[osp(2|n)]$, which
establishes our notational conventions, and an explanation of the
quasi-triangular Hopf superalgebra structure.  In Section 3 we choose an
ansatz for the Lax operator, and then determine the explicit form 
by imposing that the co-product is intertwined.  We show that the result is
consistent not only with the defining relations for the universal
$R$-matrix when the vector representation is applied to the left hand
side, but also the defining relations of the quantum superalgebra, 
including the higher-order $q$-Serre relations.  
Lastly, in Section 4 we use the Lax operator to
find the $R$-matrix for the vector representation, $\mathfrak{R} = (\pi 
\otimes \pi) \mathcal{R}$.

\section{The Construction of $U_q[osp(2|n)]$}

\noindent A quantum superalgebra is a generalisation of a classical
superalgebra involving a complex parameter $q$, which reduces to the
classical case as $q \rightarrow 1$.  In particular, $U_q[osp(2|n)]$
is a \textit{$q$-deformation} of the original enveloping algebra of
$osp(2|n)$.  Throughout this paper $q$ is assumed to be real.

We construct the quantised orthosymplectic superalgebra
$U_q[osp(2|n)]$ in a very similar way to that used for $U_q[osp(m|n)]$
in \cite{us}, but provide some details here in order to establish our
notational conventions.  For more information refer to \cite{us}.

The grading of $a$ is denoted $[a]$, where $[a] \in \mathbb{Z}_2$ is given by

\begin{equation*}
[a] =
   \begin{cases}
     0, \qquad \; a=i, & 1 \leq i \leq m, \\ 1, \qquad \; a = \mu, & 1
     \leq \mu \leq n.
   \end {cases}
\end{equation*}

\noindent Throughout this paper we use Greek letters $\mu, \nu$
etc. to denote odd indices and Latin letters $i, j$ etc. for even
indices.  If the grading is unknown, the usual $a, b, c$ etc. are
used.  Which convention applies will be clear from the context.  Throughout 
the paper we also use the symbols $\overline{a}$
and $\xi_a$, which are given by:

\begin{equation*}
\overline{a}=
  \begin{cases}
     m + 1 - a, & [a]=0, \\ n + 1 - a, & [a]=1,
  \end {cases}
\qquad \text{and} \quad \xi_{a} =
  \begin{cases}
     1, & [a] = 0, \\ (-1)^a, & [a]=1.
  \end{cases}
\end{equation*}

As a weight system, we take the set $\{ \varepsilon_1,\varepsilon_2 \}
\cup \{ \delta_\mu,\; 1 \leq \mu \leq n \}$, where
\mbox{$\varepsilon_1 = - \varepsilon_2$} and
\mbox{$\delta_{\overline{\mu}} = - \delta_\mu$}.  Acting on these
weights, we have the invariant bilinear form defined by:

\begin{equation*}
(\varepsilon_1, \varepsilon_1) = 1, \quad (\delta_\mu, \delta_\nu) =
-\delta^\mu_\nu, \quad (\varepsilon_1, \delta_\mu) = 0, \qquad 1 \leq
\mu, \nu \leq k.
\end{equation*}

\noindent When describing an object with unknown grading indexed by
$a$ the weight will be described generically as $\varepsilon_a$.  This
should not be assumed to be an even weight.

The even positive roots of $osp(2|n)$ are those of $sp(n)$, namely:

\begin{alignat*}{3}
&\delta_\mu + \delta_\nu, \quad && 1 \leq \mu,\nu \leq k, \\
&\delta_\mu - \delta_\nu, && 1 \leq \mu < \nu \leq k.
\end{alignat*}

\noindent The root system also contains a set of odd positive roots,
which are:

\begin{equation*}
\delta_\mu \pm \varepsilon_1, \qquad 1 \leq \mu \leq k.
\hspace{2cm}
\end{equation*}

\noindent Throughout this paper we choose to use the following set of
simple roots:

\begin{align*}
&\alpha_\mu = \delta_\mu - \delta_{\mu+1}, \hspace{9mm} 1 \leq \mu <
k,\notag\\ &\alpha_s = \delta_k - \varepsilon_1, \\ &\alpha_t =
\delta_k + \varepsilon_1.
\end{align*}

In $U_q[osp(m|n)]$ the graded commutator is realised by

\begin{equation*}
[A,B] = AB - (-1)^{[A][B]} BA
\end{equation*}
 
\noindent and tensor product multiplication is given by

\begin{equation*}
(A \otimes B) (C \otimes D) = (-1)^{[B][C]} (AC \otimes BD).
\end{equation*}

\noindent Using these conventions, we have:

\begin{definition} \label{def} The quantum superalgebra $U_q[osp(2|n)]$ is 
generated by simple generators $e_a, f_a, h_a$ subject to the relations:

\begin{alignat}{2}
&[h_a, e_b] = (\alpha_a, \alpha_b) e_b, && \notag \\ 
&[h_a, f_b] = - (\alpha_a, \alpha_b) f_b, && \notag \\ 
&[h_a, h_b] = 0,&& \notag \\
&[e_a, f_b] = \delta^a_b \frac{(q^{h_a} - q^{-h_a})}{(q - q^{-1})},&&
\notag \\ 
&[e_a, e_a] = [f_a,f_a]=0 & \quad &\text{for }\; (\alpha_a,\alpha_a)=0,\notag\\
&(\text{ad}\,e_b\, \circ)^{1-a_{bc}} e_c = 0 &&
  \text{for }b \neq c, \; (\alpha_b, \alpha_b)\neq 0, \label{qS1} \\
&(\text{ad}\,f_b\, \circ)^{1-a_{bc}} f_c = 0 && \text{for } b \neq c , \;
  (\alpha_b, \alpha_b) \neq 0 \label{qS2}.
\end{alignat}

\noindent Here $a_{bc}$ are the elements of the corresponding Cartan
matrix, namely:

\begin{equation*}
a_{bc} =
\begin{cases}
   \frac{2(\alpha_b, \alpha_c)}{(\alpha_b, \alpha_b)}, \quad &
(\alpha_b, \alpha_b) \neq 0, \\ (\alpha_b, \alpha_c), &(\alpha_b,
\alpha_b) = 0,
\end{cases}
\end{equation*}

\noindent and ad represents the adjoint action, which is defined
later in this section.  The relations \eqref{qS1} and \eqref{qS2} are
called the \textit{$q$-Serre relations}.  There are also extra
$q$-Serre relations which are not included here, a complete list of
which can be found in \cite{Yamane}.

We remark that $U_q[osp(m|n)]$ has the structure of a quasi-triangular
Hopf superalgebra. In 
particular, there is a linear mapping known as the \textit{coproduct},
$\Delta: U_q[osp(m|n)] \rightarrow U_q[osp(m|n)]^{\otimes
2}$, which is defined on the simple generators by:

\begin{align*}
&\Delta (e_a) = q^{\frac{1}{2}h_a} \otimes e_a + e_a \otimes
  q^{-\frac{1}{2} h_a}, \notag \\ &\Delta (f_a) = q^{\frac{1}{2}h_a}
  \otimes f_a + f_a \otimes q^{-\frac{1}{2} h_a}, \\ 
&\Delta (q^{\pm \frac{1}{2}h_a}) = q^{\pm \frac{1}{2}h_a} \otimes q^{\pm
  \frac{1}{2}h_a}, 
\end{align*}

\noindent and extends to arbitrary elements according to the homomorphism 
 property, namely:
\begin{align*}
\Delta (AB) = \Delta(A) \Delta(B). \label{coprod2}
\end{align*}

\noindent Similarly, there is an \textit{antipode}, $S: U_q[osp(2|n)]
\rightarrow U_q[osp(2|n)]$, defined on the simple generators by:

\begin{align*}
&S (e_a) = - q^{-\frac{1}{2}(\alpha_a, \alpha_a)} e_a, \notag \\ 
&S(f_a) = - q^{\frac{1}{2}(\alpha_a, \alpha_a)} f_a, \notag \\ 
&S(q^{\pm h_a}) = q^{\mp h_a}, \notag \\ 
\end{align*}

\noindent and extending to arbitrary elements according to the antihomorphism 
property, namely:

\begin{equation*}
S (ab) = (-1)^{[a][b]} S(b) S(a).
\end{equation*}
\end{definition}

It can be shown that both the coproduct and antipode are consistent
with the defining relations of the superalgebra.  These mappings are
necessary to define the adjoint action for a quantum superalgebra.  If
we adopt Sweedler's notation for the coproduct,

\begin{equation*}
\Delta(a) = \sum_{(a)} a^{(1)} \otimes a^{(2)},
\end{equation*}

\noindent the \textit{adjoint action} of $a$ on $b$ is defined to be

\begin{equation} \label{adj}
\text{ad}\; a \circ b = \sum_{(a)} (-1)^{[b][a^{(2)}]} a^{(1)} b S(a^{(2)}).
\end{equation}

One quantity that repeatedly arises in calculations for both classical
and quantum Lie superalgebras is $\rho$, the \textit{graded half-sum
of positive roots}.  In the case of $U_q[osp(2|n)]$ it is given by:

\begin{equation*}
\rho = \frac{1}{2} \sum_{\mu=1}^k (n-2\mu) \delta_\mu = \sum_{\mu=1}^k
  (k-\mu) \delta_\mu.
\end{equation*}

\noindent This satisfies the property $(\rho, \alpha) = \frac{1}{2}
(\alpha, \alpha)$ for all simple roots $\alpha$.

\subsection{$U_q[osp(2|n)]$ as a Quasi-Triangular Hopf Superalgebra}

\noindent As quantum $U_q[osp(2|n)]$ is a quasi-triangular Hopf
superalgebra, it must contain a universal $R$-matrix, providing a
solution to the quantum Yang--Baxter equation.  Before elaborating, we
need to introduce the graded twist map.

 The \textit{graded twist map} $T:U_q[osp(2|n)]^{\otimes 2}
 \rightarrow U_q[osp(2|n)]^{\otimes 2}$ is given by

\begin{equation*}
T(A \otimes B) = (-1)^{[A][B]} (B \otimes A).
\end{equation*}

\noindent For convenience $T \circ \Delta$, the twist map composed
with the coproduct, is denoted $\Delta^T$.  Then a \textit{universal
$R$-matrix}, $\mathcal{R}$, is an even, non-singular element of
$U_q[osp(2|n)]^{\otimes 2}$ satisfying the following properties:

\begin{align}
&\mathcal{R} \Delta (a) = \Delta^T (a)\mathcal{R}, \quad \forall a \in
  U_q[osp(2|n)], \notag \\ &(\text{id} \otimes \Delta) \mathcal{R} =
  \mathcal{R}_{13} \mathcal{R}_{12}, \notag \\ &(\Delta \otimes
  \text{id}) \mathcal{R} = \mathcal{R}_{13} \mathcal{R}_{23}.
  \label{Requations}
\end{align}

\noindent Here $\mathcal{R}_{yz}$ represents a copy of $\mathcal{R}$
acting on the $y$ and $z$ components respectively of $U \otimes U
\otimes U$, where each $U$ is a copy of the quantum superalgebra
$U_q[osp(2|n)]$.  When $y>z$ the usual grading term from the twist map
is included, so for example $\mathcal{R}_{21} = [\mathcal{R}^T]_{12}
$, where $\mathcal{R}^T = T (\mathcal{R})$ is the \textit{opposite
universal $R$-matrix}.

 The reason $R$-matrices are significant is that as a consequence of
 \eqref{Requations} they satisfy the quantum Yang--Baxter Equation,
 which is prominent in the study of integrable systems \cite{Baxter}:

\begin{equation*}
\mathcal{R}_{12} \mathcal{R}_{13} \mathcal{R}_{23} = \mathcal{R}_{23}
  \mathcal{R}_{13} \mathcal{R}_{12}
\end{equation*}

A superalgebra may contain many different universal $R$-matrices, but
there is a unique one belonging to $U_q[osp(2|n)]^- \otimes
U_q[osp(2|n)]^+$, and its opposite $R$-matrix in $U_q[osp(2|n)]^+
\otimes U_q[osp(2|n)]^-$.  Here $U_q[osp(2|n)]^-$ is the Hopf
subsuperalgebra generated by the lowering generators and Cartan
elements, while $U_q[osp(2|n)]^+$ is generated by the raising
generators and Cartan elements.  These particular $R$-matrices arise
out of the $\mathbb{Z}_2$-graded version of Drinfeld's double construction 
\cite{gzb}.  In this paper we consider the universal $R$-matrix belonging 
to $U_q[osp(2|n)]^-\otimes U_q[osp(2|n)]^+$.

\section{An Ansatz for the Lax Operator}

\noindent Now we have the information necessary to construct a
\textit{Lax operator} for $U_q[osp(2|n)]$.  Previously this has been
done in the superalgebra case for $U_q[gl(m|n)]$ \cite{Zhang2} and
$U_q[osp(m|n)],\, m \geq 2$ \cite{us}. First we need to define the
vector representation.

Let $\text{End} \; V$ be the space of endomorphisms of
$V$, an $(m+n)$-dimensional vector space.  Then the irreducible
\textit{vector representation} $\pi: U_q[osp(m|n)] \rightarrow
\text{End} \; V$ acts on the $U_q[osp(m|n)]$ generators as 
given in Table \ref{vectorrep}, where $E^a_b$ is the elementary matrix with a 
1 in the $(a,b)$ position and zeroes elsewhere.

\begin{table}[ht]
\caption{The action of the vector representation $\pi$ on the simple 
generators of $U_q[osp(m|n)]$} 
\label{vectorrep}
\centering
\begin{tabular}{|l|l|l|l|}\hline
\multicolumn{1}{|c|}{$\alpha_a$}& \multicolumn{1}{c|}{$\pi(e_a)$} & 
  \multicolumn{1}{c|}{$\pi(f_a)$} & \multicolumn{1}{c|}{$\pi(h_a)$} \\ \hline 
$\alpha_\mu, 1 \leq \mu <k$ & 
  $E^\mu_{\mu+1} + E^{\overline{\mu+1}}_{\overline{\mu}}$&
  $E_\mu^{\mu+1} + E_{\overline{\mu+1}}^{\overline{\mu}}$&
  $E^{\mu+1}_{\mu+1} - E^{\overline{\mu+1}}_{\overline{\mu+1}}
               - E^\mu_\mu + E^{\overline{\mu}}_{\overline{\mu}}$ \\
$\alpha_s$ &
  $E^{\mu=k}_{i=1} + (-1)^k E^{\overline{i=1}}_{\overline{\mu=k}}$&
  $- E^{i=1}_{\mu=k} + (-1)^k E_{\overline{i=1}}^{\overline{\mu=k}}$&
  $- E^{i=1}_{i=1} + E^{\overline{i}=\overline{1}}_{\overline{i}=\overline{1}} 
   -  E^{\mu=k}_{\mu=k} + E^{\overline{\mu}=\overline{k}}_
           {\overline{\mu}=\overline{k}}$ \\ 
$\alpha_t$ & $E^{\mu=k}_{\overline{i=1}} + (-1)^k E^{i=1}_{\overline{\mu=k}}$&
  $- E^{\overline{i=1}}_{\mu=k} + (-1)^k E_{i=1}^{\overline{\mu=k}}$&
  $ E^{i=1}_{i=1} - E^{\overline{i}=\overline{1}}_{\overline{i}=\overline{1}} 
   -  E^{\mu=k}_{\mu=k} + E^{\overline{\mu}=\overline{k}}_
           {\overline{\mu}=\overline{k}}$ \\ \hline
\end{tabular}
\end{table}

Now let $\mathcal{R}$ be a universal $R$-matrix of $U_{q}[osp(m|n)]$. 
 The Lax operator associated with $\mathcal{R}$ is given by

\begin{equation*}
R = (\pi \otimes \text{id}) \mathcal{R} \in (\text{End} \; V) \otimes
U_{q}[osp(m|n)]
\end{equation*}

\noindent and the $R$-matrix in the vector representation
$\mathfrak{R}$ is given by:

\begin{equation*}
\mathfrak{R} = (\pi \otimes \pi) \mathcal{R} = (\text{id} \otimes \pi)
R \in (\text{End} \; V) \otimes (\text{End} \; V).
\end{equation*}

\noindent Then the Yang-Baxter equation reduces to:

\begin{equation*}
\mathfrak{R}_{12} R_{13} R_{23} = R_{23} R_{13} \mathfrak{R}_{12}
\end{equation*}

\noindent acting on the space $V \otimes V \otimes U_q[osp(m|n)]$.

In the following sections we make use of the bra and ket notation.
The set $\{| a \rangle, a=1,...,n+2 \}$ is a basis for $V$ satisfying
the property

\begin{equation*}
E^a_b |c \rangle = \delta^c_b |a \rangle.
\end{equation*}

\noindent The set $\{ \langle a|,a=1,...,n+2 \}$ is the dual basis
such that

\begin{equation*}
\langle c| E^a_b = \delta^a_c \langle b| \quad \text{and} \quad
\langle a|b \rangle = \delta^a_b.
\end{equation*}   

As we wish to find the Lax operator belonging to $ \pi \bigl(
U_q[osp(2|n)]^-\bigr) \otimes U_q[osp(2|n)]^+$, we adopt the following
ansatz for $R$:

\begin{equation*} \label{RmatrixAnsatz}
R \equiv q ^ {\underset{a}{\sum} \pi(\tilde{h}_{a}) \otimes \tilde{h}^{a}} 
\Bigl[ I \otimes I + (q - q^{-1}) \sum_{\varepsilon_{a} < \varepsilon_{b}}
(-1)^{[b]} E^a_b \otimes \hat{\sigma}_{ba} \Bigr].
\end{equation*} 

\noindent Here $\{ \tilde{h}_a \}$ is a basis for the Cartan subalgebra and 
$\{ \tilde{h}^a \}$ the dual basis. The $\hat{\sigma}_{ba}$ are
unknown operators for which we are trying to solve.  Throughout this
paper, when working in the vector representation, we simply use $h_a$
rather than $\pi(h_a)$, and $e_a$ rather than $\pi(e_a)$.
\subsection{Constraints Arising from the Defining Relations} 
  \label{Developing Relations}

\noindent The Lax operator $R$ must be consistent with the defining
relations for the $R$-matrix, which were given as equation
\eqref{Requations}.  In particular, it must satisfy the intertwining
property for the raising generators,

\begin{equation*}
R \Delta (e_c) = \Delta^T (e_c) R.
\end{equation*}

\noindent To apply this, recall that

\begin{equation*}
\Delta(e_c) = q^{\frac{1}{2} h_c} \otimes e_c + e_c \otimes
  q^{-\frac{1}{2} h_c}.
\end{equation*}

\noindent Then, from the defining relations, we have

\begin{align*}
\Delta^T(e_c) q^{\underset{a}{\sum} \tilde{h}_a \otimes \tilde{h}^a} &=
  (q^{\frac{1}{2} h_c} \otimes e_c + e_c \otimes q^{-\frac{1}{2} h_c})
  q^{\underset{a}{\sum} \tilde{h}_a \otimes \tilde{h}^a} \\ 
&= q^{\underset{a}{\sum} \tilde{h}_a \otimes \tilde{h}^a} (e_c \otimes 
  q^{-\frac{1}{2} h_c} + q^{-\frac{3}{2} h_c} \otimes e_c).
\end{align*}

\noindent Using this, we see

\begin{align}
\Delta^T(e_c)R &= q^{\underset{a}{\sum} \tilde{h}_a \otimes \tilde{h}^a} (e_c 
  \otimes q^{-\frac {1}{2} h_c} + q^{-\frac{3}{2} h_c} \otimes e_c) 
  \Bigl[ I \otimes I + (q - q^{-1})
  \sum_{\varepsilon_{a}< \varepsilon_{b}} (-1)^{[b]} E^a_b \otimes
  \hat{\sigma}_{ba} \Bigr] \notag \\ &= q^{\underset{a}{\sum} \tilde{h}_a
  \otimes \tilde{h}^a} \biggl\{ e_c \otimes q^{-\frac{1}{2} h_c} +
  q^{-\frac{3}{2} h_c} \otimes e_c \notag \\ 
&\quad + (q-q^{-1}) \sum_{\varepsilon_{a} < \varepsilon_{b}} (-1)^{[b]}
  \Bigl[e_c E^a_b \otimes q^{-\frac{1}{2} h_c} \hat{\sigma}_{ba} 
  + (-1)^{([a] + [b])[c]} q^{-\frac{3}{2} (\alpha_c,
  \varepsilon_a)} E^a_b \otimes e_c \hat{\sigma}_{ba} \Bigr] \biggr\}.
   \label{delR} 
\end{align}

\noindent Also,

\begin{align}
R \Delta (e_c) &= q^{\underset{a}{\sum}\tilde{h}_a \otimes \tilde{h}^a} 
  \biggl\{q^{\frac{1}{2} h_c} \otimes e_c + e_c \otimes q^{-\frac{1}{2} h_c}
   \notag \\ & \hspace{2cm} + (q - q^{-1}) \sum_{\varepsilon_{a} <
   \varepsilon_{b}} (-1)^{[b]} \Bigl[ q^{\frac{1}{2} (\alpha_c,
   \varepsilon_b)} E^a_b \otimes \hat{\sigma}_{ba} e_c 
   + (-1)^{([a]+[b])[c]} E^a_b e_c \otimes
   \hat{\sigma}_{ba} q^{-\frac{1}{2} h_c} \Bigr] \biggr\}
   \label{Rdel}.
\end{align}


Hence to apply the intertwining property we simply equate (\ref{delR})
and (\ref{Rdel}).  First note that $R$ is weightless, so
$\hat{\sigma}_{ba}$ has weight $\varepsilon_b - \varepsilon_{a}$, and
thus

\begin{equation*}
q^{-\frac{1}{2} h_c} \hat{\sigma}_{ba} = q^{-\frac{1}{2} (\alpha_c,
  \varepsilon_b - \varepsilon_a)} \hat{\sigma}_{ba} q^{-\frac{1}{2}
  h_c}.
\end{equation*}

\noindent Then, equating those terms with zero weight in the first
element of the tensor product, we obtain

\begin{equation}
(q^{\frac{1}{2} h_c} - q^{-\frac{3}{2} h_c}) \otimes e_c =
  (q-q^{-1}) \sum_{\varepsilon_{b} - \varepsilon_{a} = \alpha_c}
  (-1)^{[b]} \bigl( q^{-\frac{1}{2} (\alpha_c, \alpha_c)} e_c E^a_b -
  (-1)^{[c]} E^a_b e_c \bigr) \otimes \hat{\sigma}_{ba}
  q^{-\frac{1}{2} h_c}.\label{eq1}
\end{equation}

\noindent Comparing the remaining terms, we also find

\begin{multline} \label{eq2}
\underset{\varepsilon_{b} - \varepsilon_{a} \neq \alpha_c}
  {\sum_{\varepsilon_ {a} < \varepsilon_{b}}} (-1)^{[b]} \bigl(
  q^{-\frac{1}{2}(\alpha_c, \varepsilon_{b} - \varepsilon_{a})} e_c
  E^a_b - (-1)^{([a] + [b])([c])} E^a_b e_c \bigr) \otimes
  \hat{\sigma}_{ba} q^{-\frac{1}{2} h_c} \\ = \sum_{\varepsilon_{a} <
  \varepsilon_{b}} (-1)^{[b]} E^a_b \otimes \bigl( q^{\frac{1}{2}
  (\alpha_c, \varepsilon_b)} \hat{\sigma}_{ba} e_c -
  (-1)^{([a]+[b])[c]} q^{-\frac{3}{2} (\alpha_c, \varepsilon_a)} e_c
  \hat{\sigma}_{ba} \bigr).
\end{multline}



From the first of these equations we can deduce the simple values of
$\hat{\sigma}_{ba}$; from the second, relations involving all the
$\hat{\sigma}_{ba}$.  Before doing so, however, it is convenient to
define a new set, $\overline{\Phi}^+$.

\begin{definition}
The extended system of positive roots, $\overline{\Phi}^+$, is defined
by

\begin{equation*}
\overline{\Phi}^+ \equiv \{\varepsilon_b - \varepsilon_a|
\varepsilon_b > \varepsilon_a\} = \Phi^+ \cup \{2 \varepsilon_1\}
\end{equation*}

\noindent where $\Phi^+$ is the usual system of positive roots.
\end{definition}

Now consider equation (\ref{eq2}).  In the case when $\varepsilon_b -
\varepsilon_a + \alpha_c \notin \overline{\Phi}^+$, by collecting the
terms of weight $\varepsilon_b - \varepsilon_a + \alpha_c$ in the
second half of the tensor product we find:

\begin{equation} \label{*}
q^{\frac{1}{2}(\alpha_c, \varepsilon_{b})} \hat{\sigma}_{ba} e_c -
    (-1)^{([a]+[b])[c]} q^{-\frac{3}{2} (\alpha_c, \varepsilon_a)} e_c
    \hat{\sigma}_{ba} = 0.
\end{equation}

\noindent Similarly, when $\varepsilon_b > \varepsilon_a$ and
$\varepsilon_b - \varepsilon_a + \alpha_c = \varepsilon_{b'} -
\varepsilon_{a'} \in \overline{\Phi}^+$ we find:

\begin{multline*}
\underset{\varepsilon_{b} - \varepsilon_{a} + \alpha_c =
  \varepsilon_{b'} - \varepsilon_{a'}} {\sum_{\varepsilon_{a'} <
  \varepsilon_{b'}}} \hspace{-8mm} (-1)^{[b']}
  \bigl(q^{-\frac{1}{2}(\alpha_c, \varepsilon_{b'}-\varepsilon_{a'})}
  e_c E^{a'}_{b'} - (-1)^{([a']+[b'])[c]} E^{a'}_{b'} e_c \bigr)
  \otimes \hat{\sigma}_{b'\!a'} q^{-\frac{1}{2} h_c} \\ = (-1)^{[b]}
  E^a_b \otimes \bigl( q^{\frac{1}{2} (\alpha_c, \varepsilon_b)}
  \hat{\sigma}_{ba} e_c - (-1)^{([a]+[b])[c])} q^{-\frac{3}{2}
  (\alpha_c, \varepsilon_a)} e_c \hat{\sigma}_{ba} \bigr).
\end{multline*}

\noindent However $e_c E^{a'}_{b'}$ and $E^a_{b}$ are linearly
independent unless $b=b'$, as are $E^{a'}_{b'} e_c$ and $E^{a}_{b}$
for $a \neq a'$, and thus this equation reduces to

\begin{align*}
q^{-\frac{1}{2} (\alpha_c, \varepsilon_b - \varepsilon_{a} +
   \alpha_c)}& e_c E^{a'}_{b} \otimes \hat{\sigma}_{ba'}
   q^{-\frac{1}{2} h_c} \Big\vert_{\varepsilon_{a'} = \varepsilon_a -
   \alpha_c} - (-1)^{([a]+[b])[c]} E^{a}_{b'}e_c \otimes
  \hat{\sigma}_{b'a} q^{-\frac{1}{2}h_c} \Big\vert_{\varepsilon_{b'} =
  \varepsilon_b +\alpha_c}\\ 
&= E^{a}_{b} \otimes \bigl( q^{\frac{1}{2} (\alpha_c,\varepsilon_b)} 
  \hat{\sigma}_{ba} e_c - (-1)^{([a]+[b])[c]} q^{-\frac{3}{2} 
  (\alpha_c,\varepsilon_a)} e_c \hat{\sigma}_{ba} \bigr)
\end{align*}

\noindent for $\varepsilon_b > \varepsilon_a$.  This further
simplifies to

\begin{multline}
q^{-\frac{1}{2} (\alpha_c, \alpha_c - \varepsilon_{a})} \langle
  a|e_c|a' \rangle \hat{\sigma}_{ba'}-(-1)^{([a]+[b])[c]}
  q^{\frac{1}{2} (\alpha_c, \varepsilon_b)} \langle b'|e_c|b \rangle
  \hat{\sigma}_{b'a} \\ = q^{(\alpha_c,\varepsilon_b)}
  \hat{\sigma}_{ba}e_c q^{\frac{1}{2} h_c} - (-1)^ {([a]+[b])[c]}
  q^{-(\alpha_c,\varepsilon_a)} e_c q^{\frac{1}{2} h_c}
  \hat{\sigma}_{ba} \label{**}
\end{multline}

\noindent for $\varepsilon_b > \varepsilon_a$.  All the necessary
information is contained within equations \eqref{eq1} and \eqref{**}.
To construct the Lax operator $R = (\pi \otimes 1) \mathcal{R}$ first
we use equation \eqref{eq1} to find the solutions for
$\hat{\sigma}_{ba}$ associated with the simple roots $\alpha_c$.  Then
we apply the recursion relations arising from \eqref{**} to find the
remaining values of $\hat{\sigma}_{ba}$.


\subsection{The Simple Operators} \label{iv}

\noindent In this section we solve equation \eqref{eq1} to find the simple 
values of $\hat{\sigma}_{ba}$, namely those for which $\varepsilon_b - 
\varepsilon_a$ is a simple root.  To solve the equation we must consider the 
various simple roots individually.  Consider the odd root $\alpha_i = \delta_k
 +\varepsilon_1$.  In the vector representation $e_t =
 E^{\mu=k}_{\overline{i}=\overline{1}} +(-1)^k
 E^{i=1}_{\overline{\mu}=\overline{k}}$ and $h_t = E^{i=1}_{i=1} -
 E^{\overline{i}=\overline{1}}_{\overline{i}=\overline{1}} -
 E^{\mu=k}_{\mu=k} +
 E^{\overline{\mu}=\overline{k}}_{\overline{\mu}=\overline{k}}$.
 Hence the left-hand side of (\ref{eq1}) becomes:

\begin{align*}
LHS &= (q^{\frac{1}{2} h_t} - q^{-\frac{3}{2} h_t}) \otimes e_t \\ &=
    (q-q^{-1}) \bigl\{ q^{-\frac{1}{2}} (E^{i=1}_{i=1} +
    E^{\overline{\mu}=\overline{k}}_{\overline{\mu}=\overline{k}}) -
    q^{\frac{1}{2}}
    (E^{\overline{i}=\overline{1}}_{\overline{i}=\overline{1}} +
    E^{\mu=k}_{\mu=k}) \bigr\} \otimes e_t,
\end{align*}

\noindent whereas the right-hand side is:

\begin{align*}
RHS &= (q-q^{-1}) \sum_{\varepsilon_{b} - \varepsilon_{a} = \alpha_t}
   (-1)^{[b]} \bigl( e_t E^a_b + E^a_b e_t \bigr) \otimes
   \hat{\sigma}_{ba} q^{-\frac{1}{2} h_t} \notag \\ &= (q-q^{-1})
   \Bigl\{ (-1)^{k}( E^{i=1}_{i=1} + E^{\overline{\mu}=\overline{k}}
   _{\overline{\mu}=\overline{k}}) \otimes
   \hat{\sigma}_{(i=1)(\overline{\mu}=\overline{k})} q^{-\frac{1}{2} h_t} 
   (E^{\overline{i}=\overline{1}}_{\overline{i}=\overline{1}} +
   E^{\mu=k}_{\mu=k}) \otimes \hat{\sigma}_{(\mu=k) (\overline{i}=
   \overline{1})} q^{-\frac{1}{2} h_t} \Bigr\}.
\end{align*}

\noindent Equating these gives

\begin{equation*}
\hat{\sigma}_{(\mu=k)(\overline{i}=\overline{1})}= (-1)^k q
  \hat{\sigma}_{(i=1)(\overline{\mu}=\overline{k})} = q^{\frac{1}{2}}
  e_t q^{\frac{1}{2} h_t}.
\end{equation*}

 By performing similar calculations for the other simple roots we
 obtain the simple operators given below in Table \ref{fundval}.
 These values for $\hat{\sigma}_{ba}$ form the basis for finding $R$,
 as from these all the others can be explicitly determined in any
 given representation.

\begin{table}[ht]
\caption{The simple values for $\hat{\sigma}_{ba}$.} \label{fundval}
\centering
\begin{tabular}{|l|lll|} \hline
\multicolumn{1}{|c|}{Simple Root}& \multicolumn{3}{c|}{Corresponding
  $\hat{\sigma}_{ba}$} \\ \hline $\alpha_\mu = \delta_\mu -
  \delta_{\mu+1},\,\mu < k$ & $\hat{\sigma}_{\mu (\mu+1)} $&$ =
  \hat{\sigma}_{(\overline{\mu+1}) \overline{\mu}} $&$=
  q^{-\frac{1}{2}} e_\mu q^{\frac{1}{2} h_\mu}$ \\ $\alpha_s =
  \delta_k - \varepsilon_1$ & $\hat{\sigma}_{(\mu=k) (i=1)} $&$ = (-1)^k
  q \hat{\sigma}_{(\overline{i}=\overline{1})(\overline{\mu} =
  \overline{k})} $&$ = q^{\frac{1}{2}} e_s q^{\frac{1}{2} h_s}$\\
  $\alpha_t = \delta_k + \varepsilon_1$ &$\hat{\sigma}_{(\mu=k0
  (\overline{i}= \overline{1})} $&$ = (-1)^k q
  \hat{\sigma}_{(i=1)(\overline{\mu}=\overline{k})} $&$ =
  q^{\frac{1}{2}} e_t q^{\frac{1}{2} h_t}$ \\ \hline
\end{tabular}
\end{table}

\subsection{Constructing the Non-Simple Operators} \label{noniv}

\noindent Now we develop the recurrence relations required to
calculate the remaining values of $\hat{\sigma}_{ba}$.  Consider equation 
\eqref{**}. To extract the recurrence relations to be applied to the
 simple values of $\hat{\sigma}_{ba}$, we must again consider the
 simple roots individually.  We consider the case $\alpha_t = \delta_k
 + \varepsilon_1$, so $e_t =
 \sigma^{\mu=k}_{\overline{i}=\overline{1}} \equiv E^{\mu=k}
 _{\overline{i}=\overline{1}} +(-1)^k
 E^{i=1}_{\overline{\mu}=\overline{k}}$.  Now

\begin{equation*}
\langle a|e_t = \delta^a_{\mu=k} \langle \overline{i}=\overline{1}| -
  \delta^a_{i=1} \langle \overline{\mu}=\overline{k}|,\qquad e_t |b
  \rangle = \delta^b_{\overline{i}=\overline{1} } |\mu=k \rangle -
  \delta^b_{\overline{\mu}=\overline{k}} |i=1 \rangle.
\end{equation*}

\noindent We then apply this to equation \eqref{**} to obtain:

\begin{align*}
q^{-\frac{1}{2}} \delta^a_{\mu=k} \hat{\sigma}_{b(\overline{i} =
   \overline{1})} &+ (-1)^k q^{\frac{1}{2}} \delta^a_{i=1}
   \hat{\sigma}_{b(\overline{\mu}=\overline{k})} -(-1)^{[a]}
   q^{-\frac{1}{2}} \delta^b_{\overline{i} = \overline{1}}
   \hat{\sigma}_{(\mu=k) a} + (-1)^{[a]} (-1)^k q^{\frac{1}{2}}
   \delta^b_{\overline{\mu} = \overline{k}} \hat{\sigma}_{(i=1) a} \\ 
&= q^{(\alpha_t,\varepsilon_b)} \hat{\sigma}_{ba} e_t
   q^{\frac{1}{2} h_t} - q^{-(\alpha_t,\varepsilon_a) } e_t
   q^{\frac{1}{2} h_t} \hat{\sigma}_{ba}
\end{align*}

\noindent for $\varepsilon_b >\varepsilon_a.$ Recalling that
 $\hat{\sigma}_{(\mu=k)(\overline{i}=\overline{1})}= (-1)^k q
 \hat{\sigma}_{(i=1)(\overline{\mu}=\overline{k})} = q^{\frac{1}{2}} e_t
 q^{\frac{1}{2} h_t} $, we can deduce the following relations:

\begin{alignat}{2}
&\hat{\sigma}_{\nu (\overline{i}=\overline{1})} =
   \hat{\sigma}_{\nu(\mu=k)} \hat{\sigma}_{(\mu=k)(\overline{i}=\overline{1})} 
 - q \hat{\sigma}_{(\mu=k)(\overline{i}=\overline{1})} 
   \hat{\sigma}_{\nu(\mu=k)} , & \qquad &\nu<k, \notag \\ 
&\hat{\sigma}_{(\mu=k)\overline{\nu}} = q^{(\delta_k, \delta_\nu)} 
   \hat{\sigma}_{(\mu=k)(\overline{i}=\overline{1})}  
    \hat{\sigma}_{(\overline{i}=\overline{1})\overline{\nu}} + q^{-1}
   \hat{\sigma}_{(\overline{i}=\overline{1})\overline{\nu}}
   \hat{\sigma}_{(\mu=k)(\overline{i}=\overline{1})}, & &\nu \leq k,\notag \\ 
&\hat{\sigma}_{(i=1)\, \overline{\nu}} 
   = \hat{\sigma}_{(i=1)(\overline{\mu}= \overline{k})}
   \hat{\sigma}_{(\overline{\mu}=\overline{k})\overline{\nu}} - q
   \hat{\sigma}_{(\overline{\mu}=\overline{k})\overline{\nu}}
   \hat{\sigma}_{(i=1)(\overline{\mu}=\overline{k})}, && \nu < k \notag \\ 
&\hat{\sigma}_{\nu(\overline{\mu}=\overline{k})} =
   q^{(\delta_k, \delta_\nu)} \hat{\sigma}_{\nu (i=1)}
   \hat{\sigma}_{(i=1) (\overline{\mu}=\overline{k})} + q^{-1}
   \hat{\sigma}_{(i=1) (\overline{\mu}=\overline{k})}
   \hat{\sigma}_{\nu (i=1)}, && \nu \leq k \notag \\ 
&q^{(\alpha_t,   \varepsilon_b)} \hat{\sigma}_{ba} \hat{\sigma}_{(\mu=k)
   (\overline{i}=\overline{1})} - (-1)^{[a]+[b]} q^{-(\alpha_t,
   \varepsilon_a)} \hat{\sigma}_{(\mu=k)(\overline{i}=\overline{1})}
   \hat{\sigma}_{ba} = 0, && \varepsilon_b > \varepsilon_a;\;
   \varepsilon_a \neq \varepsilon_1, \delta_k \notag \\ &&& \text{and
   }\varepsilon_b \neq -\varepsilon_1, -\delta_k. \notag
\end{alignat}

 We then follow the same procedure to find the relations associated
 with the other simple roots, which are listed in Table \ref{list} in
 Appendix A.  The relations can be summarised in a compact form.
 There are two different types of relations; recursive and
 $q$-commutative.  The latter can be condensed into:

\vspace{-2mm}
\begin{equation} \label{commutationrelations}
\boxed{q^{(\alpha_c, \varepsilon_b)} \hat{\sigma}_{ba} e_c
  q^{\frac{1}{2} h_c} - (-1)^{([a]+[b])[c]} q^{-(\alpha_c,
    \varepsilon_a)} e_c q^{\frac{1}{2}h_c} \hat{\sigma}_{ba} = 0,
  \quad \varepsilon_b > \varepsilon_a}
\end{equation}
      
\noindent where neither $\varepsilon_a - \alpha_c$ nor $\varepsilon_b
+ \alpha_c$ equals any $\varepsilon_x$.  Note this is almost the same
as equation \eqref{*}, with slightly softer restrictions on $a$ and
$b$ the only difference.

 With one exception the recursion relations can be summarised as:

\vspace{-2mm}
\begin{equation} \label{nice}
\boxed{\hat{\sigma}_{ba} = q^{-(\varepsilon_b,\varepsilon_a)}
  \hat{\sigma}_{bc} \hat{\sigma}_{ca} - q^{-(\varepsilon_c,
  \varepsilon_c)} (-1)^{([b]+[c]) ([a]+[c])} \hat{\sigma}_{ca}
  \hat{\sigma}_{bc}, \quad \varepsilon_b > \varepsilon_c >
  \varepsilon_a}
\end{equation}

\noindent where $ c \neq \overline{b} \text{ or } \overline{a}$.  The
remaining relation is

\begin{equation} \label{2l}
\hat{\sigma}_{(\mu=k)(\overline{i} = \overline{1})} - (-1)^k q
  \hat{\sigma}_ {(i=1)(\overline{\mu} = \overline{k})} = q^{-1} [
  \hat{\sigma}_{(\mu=k)(i=2)}, \hat{\sigma}_{(i=1)(\overline{i} =
  \overline{1})}].
\end{equation}

\

\noindent Noting from Table \ref{fundval} that $\hat{\sigma}_{(\mu=k)
(\overline{i} = \overline{1})} - (-1)^k q \hat{\sigma}_{(i=1)
(\overline{\mu} = \overline{k})} = 0$ and from the Appendix that
$\hat{\sigma}_{(i=1)(\overline{i}=\overline{1})}$ commutes with all the
other simple generators, equation \eqref{2l} can be replaced with the
condition:

\begin{equation} \label{l-1lbar}
\boxed{\hat{\sigma}_{(i=1)(\overline{i}=\overline{1})} =0.}
\end{equation}

\noindent This unified form of the relations can be applied to the
simple operators given in Table \ref{fundval} to obtain all the
remaining values of $\hat{\sigma}_{ba}$ in a given representation.

 Hence we have found the following result:

\begin{proposition} \label{R1}

\noindent There is a unique operator-valued matrix $R \in (\text{End
}V) \otimes U_q[osp(2|n)]^+$ of the form

\begin{equation*}
R = q^{\underset{a}{\sum} h_{a} \otimes h^{a}} \Bigl[I \otimes I + (q
- q^{-1}) \sum_{\varepsilon_{a} < \varepsilon_{b}} (-1)^{[b]} E^a_b
\otimes \hat{\sigma}_{ba} \Bigr],
\end{equation*} 

\noindent satisfying \mbox{$R \Delta (e_c) = \Delta^T (e_c) R$}. The
simple operators for that matrix are given by:

\begin{alignat}{2}
&\hat{\sigma}_{\mu (\mu +1)} =
  \hat{\sigma}_{(\overline{\mu+1})\overline{\mu}} = q^{-\frac{1}{2}}
  e_\mu q^{\frac{1}{2} h_\mu}, &\quad &1 \leq \mu <k, \notag \\
 &\hat{\sigma}_{(\mu=k)(i=1)} = (-1)^k q\, \hat{\sigma}_{(i =\overline{1}) 
  (\overline{\mu} = \overline{k})} =q^{\frac{1}{2}} e_s q^{\frac{1}{2} h_s},\\ 
 &\hat{\sigma}_{(\mu=k)(\overline{i}=\overline{1})} = (-1)^k q
  \hat{\sigma}_{(i=1)(\overline{\mu}=\overline{k})} = q^{\frac{1}{2}}
  e_t q^{\frac{1}{2} h_t}, \\ 
 &\hat{\sigma}_{(i=1)(\overline{i}=\overline{1})} = 0;
\label{ivalues}
\end{alignat}

\noindent and the remaining values can be calculated using

\

\noindent (i) the $q$-commutation relations

\begin{equation} \label{qcom}
q^{(\alpha_c, \varepsilon_b)} \hat{\sigma}_{ba} e_c q^{\frac{1}{2}
  h_c} - (-1)^{([a]+[b])[c]} q^{-(\alpha_c, \varepsilon_a)} e_c
  q^{\frac{1}{2}h_c} \hat{\sigma}_{ba} = 0, \quad \varepsilon_b >
  \varepsilon_a,
\end{equation}
      
\noindent where neither $\varepsilon_a - \alpha_c$ nor $\varepsilon_b
+ \alpha_c$ equals any $\varepsilon_x$; and

\

\noindent (ii) the induction relations

\begin{equation} \label{indrel}
\hat{\sigma}_{ba} = q^{-(\varepsilon_b,\varepsilon_a)}
  \hat{\sigma}_{bc} \hat{\sigma}_{ca} - q^{-(\varepsilon_c,
  \varepsilon_c)} (-1)^{([b]+[c]) ([a]+[c])} \hat{\sigma}_{ca}
  \hat{\sigma}_{bc}, \quad \varepsilon_b > \varepsilon_c >
  \varepsilon_a,
\end{equation}

\noindent where $ c \neq \overline{b} \text{ or } \overline{a}$.
\end{proposition}

\

\noindent This set of simple operators and relations uniquely define
the unknowns $\hat{\sigma}_{ba}$.  To be certain that $R$ as
defined by the operators $\hat{\sigma}_{ba}$ is indeed a Lax
operator, it is necessary to check the remaining $R$-matrix properties, namely

\begin{equation} \label{1tensordel}
(\text{id} \otimes \Delta) R = R_{13} R_{12}
\end{equation}

\noindent and the intertwining property for the remaining generators,

\begin{equation*}
R \Delta(a) = \Delta^T(a) R, \qquad \forall a \in U_q[osp(2|n)].
\end{equation*}

\noindent However the calculations proceed in exactly the same manner as
those given for $m>2$ in \cite{me}, so they need not be reproduced here.  Hence
we have:

\begin{proposition}
The operator $R \in (\text{End }V) \otimes U_q[osp(2|n)]^+$ defined in
Proposition \ref{R1} is a Lax operator.
\end{proposition} 

 We now calculate the opposite Lax operator $R^T$, and briefly examine
 whether the defining relations for the $\hat{\sigma}_{ba}$
 incorporate the $q$-Serre relations for $U_q[osp(2|n)]$.

\subsection{The Opposite Lax Operator} \label{opposite}

\noindent Having found the Lax operator $R = (\pi \otimes \text{id})
\mathcal{R}$, we use that result to find its opposite $R^T = (\pi
\otimes \text{id}) \mathcal{R}^T$, where $\mathcal{R}^T$ is the
opposite universal $R$-matrix of $U_q[osp(2|n)]$.  We begin by introducing the 
graded conjugation action.

 A graded conjugation on $U_q[osp(2|n)]$, $\dagger$, is defined on the simple
 generators by:

\begin{equation*}
e_a^\dagger = f_a, \qquad f_a^\dagger = (-1)^{[a]} e_a, \qquad
h_a^\dagger = h_a.
\end{equation*}

\noindent It is consistent with the coproduct and extends naturally to
all remaining elements of $U_q[osp(2|n)]$, satisfying the properties:

\begin{align*}
& (AB)^\dagger = (-1)^{[A][B]} B^\dagger A^\dagger, \\ 
& (A \otimes B)^\dagger = A^\dagger \otimes B^\dagger, \\ 
& \Delta (A)^\dagger =\Delta (A^\dagger).
\end{align*}

 Now it was shown in \cite{us} that the opposite $R$-matrix $R^T$ is in fact 
given by $R^\dagger$. Thus we can find the opposite Lax operator $R^T$ simply 
by using the usual rules for graded conjugation.  We find the opposite Lax 
operator $R^T$ can be written as

\begin{equation} \label{RT}
R^T = \sum_a E^a_a \otimes q^{h_{\varepsilon_a}} + (q-q^{-1})
\sum_{\varepsilon_b > \varepsilon_a} (-1)^{[a]} E^b_a \otimes
\hat{\sigma}_{ab} q^{h_{\varepsilon_a}},
\end{equation}

\noindent where

\begin{equation*}
\hat{\sigma}_{ab} = (-1)^{[b]([a]+[b])} \hat{\sigma}_{ba}^\dagger,
\quad \varepsilon_b > \varepsilon_a.
\end{equation*}


\subsection{q-Serre Relations}
\noindent Having shown the relations found in Section \ref{noniv}
define a Lax operator, we also wish to see if they incorporate the
$q$-Serre relations.  It is too space-consuming to list all of these,
so we will merely provide a couple of examples, including the extra
$q$-Serre relations.

 First recall that if $\varepsilon_b - \varepsilon_a$ is a simple
 root, then $\hat{\sigma}_{ba} \propto e_c q^{\frac{1}{2}h_c}$ for
 either $c = b$ or $c=\overline{a}$.  Then setting $E_a = e_a
 q^{\frac{1}{2} h_a}$, we see from the $U_q[osp(2|n)]$ defining
 relations in Definition \ref{def} that:

\begin{alignat}{2}
&&\Delta (E_a) &= q^{h_a} \otimes E_a + E_a \otimes 1 \notag \\
&&S(E_a) &= -q^{-\frac{1}{2}(\alpha_a, \alpha_a)}
q^{-\frac{1}{2}h_a}e_a \notag \\ &&&= - q^{-h_a} E_a \notag \\
&\therefore \qquad & \text{ad}\, E_a \circ b &= - (-1)^{[a][b]} q^{h_a}b
q^{-h_a}E_a + E_a b \notag \\ &&&= E_a b - (-1)^{[a][b]}
q^{(\alpha_a,\varepsilon_b)} bE_a. \label{adjoint}
\end{alignat}

\noindent Now consider the simple generators
$\hat{\sigma}_{(\nu=k-1)(\mu=k)}$ and $\hat{\sigma}_{(\mu=k)
(\overline{i}=\overline{1})}$.

\begin{align*}
(\text{ad}\: \hat{\sigma}_{(\nu=k-1)(\mu=k)}\: \circ)^2 \hat{\sigma}_{(\mu=k)
  (\overline{i}=\overline{1})} &= \text{ad}\: \hat{\sigma}_{(\nu=k-1)(\mu=k)}
  \circ \\ &\qquad (\hat{\sigma}_{(\nu=k-1)(\mu=k)}
  \hat{\sigma}_{(\mu=k)(\overline{i}= \overline{1})} - q
  \hat{\sigma}_{(\mu=k)(\overline{i}=\overline{1})}
  \hat{\sigma}_{(\nu=k-1)(\mu=k)}) \\ 
&= \text{ad}\: \hat{\sigma}_{(\nu=k-1)(\mu=k)} \circ \hat{\sigma}_{(\nu=k-1)
  (\overline{i}=\overline{1})} \\ &= \hat{\sigma}_{(\nu=k-1)(\mu=k)}
  \hat{\sigma}_{(\nu=k-1)(\overline{i}= \overline{1})} - q^{-1}
  \hat{\sigma}_{(\nu=k-1)(\overline{i}=\overline{1})}
  \hat{\sigma}_{(\nu=k-1)(\mu=k)} \\ &= 0
\end{align*}

\noindent from \eqref{qcom}.  This is equivalent to the $q$-Serre
relation $(\text{ad}\: e_b\: \circ)^{1-a_{bc}} e_c = 0$ for this pair of
simple operators.  In a similar way, we can verify this relation for
any $b \neq c$ where $(\alpha_b,\alpha_b) \neq 0$.  The defining
relations for the $\hat{\sigma}_{ba}$, therefore, incorporate all the
standard $q$-Serre relations for raising generators.

 This still leaves the extra $q$-Serre relations, which involve the
 odd roots.  There are three of these for our choice of simple roots
 \cite{Yamane}.  Explicitly, taking into account the different
 conventions between \cite{Yamane} and here, the relevant extra
 $q$-Serre relations for $U_q[osp(2|n)]$ are:

\begin{align} 
\label{q1}
&\bigl[ \hat{\sigma}_{(\nu=k-1)(\mu=k)}, \bigl[
  \hat{\sigma}_{(\mu=k)(i=1)},
  \hat{\sigma}_{(\mu=k)(\overline{i}=\overline{1})} ]_q \, ]_q  
  \notag \\ & \hspace{35mm} 
= (q+q^{-1}) \bigl[ \hat{\sigma}_{(\mu=k)(i=1)},
  \bigl[ \hat{\sigma}_{(\nu=k-1)(\mu=k)},
  \hat{\sigma}_{(\mu=k)(\overline{i}=\overline{1})} ]_q \, ]_q \\
  &\bigl[ \hat{\sigma}_{(\nu=k-1)(\mu=k)}, \bigl[ \hat{\sigma}_{(\mu=k)
  (\overline{i}=\overline{1})}, \hat{\sigma}_{(\mu=k)(i=1)} ]_q \, ]_q
  \notag \\ & \hspace{35mm} 
= (q+q^{-1}) \bigl[
  \hat{\sigma}_{(\mu=k)(\overline{i}= \overline{1})}, \bigl[
  \hat{\sigma}_{(\nu=k-1)(\mu=k)}, \hat{\sigma}_{(\mu=k)(i=1)} ]_q \, ]_q
  \label{q2}\\ &\bigl[ \hat{\sigma}_{(\mu=k)(i=1)}, \bigl[
  \hat{\sigma}_{(\mu=k)(\overline{i}= \overline{1})},
  \hat{\sigma}_{(\nu=k-1)(\mu=k)} ]_q \, ]_q = \bigl[
  \hat{\sigma}_{(\mu=k)(\overline{i}=\overline{1})}, \bigl[
  \hat{\sigma}_{(\mu=k)(i=1)}, \hat{\sigma}_{(\nu=k-1)(\mu=k)} ]_q \, ]_q,
  \label{q3}
\end{align}

\noindent where $[x,y]_q$ represents the adjoint action $\text{ad}\, x \circ
y$.

 Consider equation \eqref{q1}.  Using the defining relations
 \eqref{qcom} and \eqref{indrel} for the $\hat{\sigma}_{ba}$ together
 with the adjoint action as given in equation \eqref{adjoint}, we
 find:

\begin{align*}
RHS &= (q+q^{-1}) \bigl[ \hat{\sigma}_{(\mu=k)(i=1)},
  \hat{\sigma}_{(\nu=k-1)(\mu=k)}
  \hat{\sigma}_{(\mu=k)(\overline{i}=\overline{1})} - q
  \hat{\sigma}_{(\mu=k)(\overline{i}=\overline{1})}
  \hat{\sigma}_{(\nu=k-1)(\mu=k)} \bigr]_q \\ 
&= (q+q^{-1}) \bigl[  \hat{\sigma}_{(\mu=k)(i=1)}, \hat{\sigma}_{(\nu=k-1)
  (\overline{i}=\overline{1})} \bigr]_q \\ 
&= (q+q^{-1}) \bigl\{  \hat{\sigma}_{(\mu=k)(i=1)} \hat{\sigma}_{(\nu=k-1)
  (\overline{i}=\overline{1})} + q^{-1} \hat{\sigma}_{(\nu=k-1)
  (\overline{i}=\overline{1})} \hat{\sigma}_{(\mu=k)(i=1)} \bigr\} \\ 
&=  (-1)^k (q+q^{-1}) \bigl\{ q \hat{\sigma}_{(\overline{i}=\overline{1}) 
  (\overline{\mu}=\overline{k})} \hat{\sigma}_{(\nu=k-1)
  (\overline{i}=\overline{1})} + \hat{\sigma}_{(\nu=k-1) 
  (\overline{i}=\overline{1})} \hat{\sigma}_{(\overline{i}=\overline{1})  
  (\overline{\mu}= \overline{k})} \bigr\}
\end{align*}

\noindent and

\begin{align*}
LHS &= \bigl[ \hat{\sigma}_{(\nu=k-1)(\mu=k)}, \hat{\sigma}_{(\mu=k)(i=1)} 
  \hat{\sigma}_{(\mu=k)(\overline{i}=\overline{1})} + q^{-2}
  \hat{\sigma}_{(\mu=k)(\overline{i}=\overline{1})}
  \hat{\sigma}_{(\mu=k)(i=1)} \bigr]_q \\ 
&=\hat{\sigma}_{(\nu=k-1)(\mu=k)} \hat{\sigma}_{(\mu=k)(i=1)}
  \hat{\sigma}_{(\mu=k)(\overline{i}=\overline{1})} - q^2
  \hat{\sigma}_{(\mu=k)(i=1)} \hat{\sigma}_{(\mu=k)(\overline{i}=\overline{1})}
  \hat{\sigma}_{(\nu=k-1)(\mu=k)} \\ 
& \quad + q^{-2} \hat{\sigma}_{(\nu=k-1)(\mu=k)} \hat{\sigma}_{(\mu=k)
  (\overline{i}=\overline{1})} \hat{\sigma}_{(\mu=k)(i=1)} -
  \hat{\sigma}_{(\mu=k)(\overline{i}=\overline{1})}
  \hat{\sigma}_{(\mu=k)(i=1)} \hat{\sigma}_{(\nu=k-1)(\mu=k)} \\ 
&= \bigl( \hat{\sigma}_{(\nu=k-1)(i=1)} + q \hat{\sigma}_{(\mu=k)(i=1)}
  \hat{\sigma}_{(\nu=k-1)(\mu=k)} \bigr) \hat{\sigma}_{(\mu=k)
  (\overline{i}=\overline{1})} \\ 
& \quad + q \hat{\sigma}_{(\mu=k)(i=1)}  \bigl( \hat{\sigma}_{(\nu=k-1)
  (\overline{i} = \overline{1})} - \hat{\sigma}_{(\nu=k-1)(\mu=k)} 
  \hat{\sigma}_{(\mu=k)(\overline{i}=\overline{1})} \bigr) \\ 
& \quad + q^{-2} \bigl( \hat{\sigma}_{(\nu=k-1)(\overline{i} = \overline{1})} 
  + q \hat{\sigma}_{(\mu=k)(\overline{i}=\overline{1})}
  \hat{\sigma}_{(\nu=k-1)(\mu=k)} \bigr) \hat{\sigma}_{(\mu=k)(i=1)} \\ 
& \quad + q^{-1} \hat{\sigma}_{(\mu=k)(\overline{i}=\overline{1})}
  \bigl( \hat{\sigma}_{(\nu=k-1)(i=1)} - \hat{\sigma}_{(\nu=k-1)(\mu=k)}
  \hat{\sigma}_{(\mu=k)(i=1)} \bigr) \\ 
&= \hat{\sigma}_{(\nu=k-1)(i=1)}  \hat{\sigma}_{(\mu=k)
  (\overline{i}=\overline{1})} + q \hat{\sigma}_{(\mu=k)(i=1)} 
  \hat{\sigma}_{(\nu=k-1)(\overline{i} = \overline{1})} \\ 
& \quad + q^{-2} \hat{\sigma}_{(\nu=k-1)(\overline{i} = \overline{1})}
  \hat{\sigma}_{(\mu=k)(i=1)} + q^{-1} \hat{\sigma}_{(\mu=k)
  (\overline{i}=\overline{1})} \hat{\sigma}_{(\mu=k)(i=1)} \\ 
&= (-1)^kq \bigl( \hat{\sigma}_{(\nu=k-1)(\overline{\mu}=\overline{k})} 
  - q^{-1} \hat{\sigma}_{(i=1)(\overline{\mu}=\overline{k})}
  \hat{\sigma}_{(\nu=k-1)(i=1)} \bigr) \\ 
& \quad + (-1)^k q^2 \hat{\sigma}_{(\overline{i}=\overline{1})(\overline{\mu} =
  \overline{k})} \hat{\sigma}_{(\nu=k-1)(\overline{i} = \overline{1})} +
  (-1)^k q^{-1} \hat{\sigma}_{(\nu=k-1)(\overline{i} = \overline{1})}
  \hat{\sigma}_{(\overline{i}=\overline{1})(\overline{\mu} =\overline{k})} \\ 
& \quad + (-1)^k \hat{\sigma}_{(i=1)(\overline{\mu}=\overline{k})} 
  \hat{\sigma}_{(\nu=k-1)(i=1)} \\ 
&= (-1)^k q \bigl( \hat{\sigma}_{(\nu=k-1)(\overline{i} = \overline{1})}
  \hat{\sigma}_{(\overline{i}=\overline{1})(\overline{\mu} = \overline{k})} 
  + q^{-1} \hat{\sigma}_{(\overline{i}=\overline{1}) 
  (\overline{\mu} = \overline{k})} \hat{\sigma}_{(\nu=k-1)
  (\overline{i} = \overline{1})} \bigr) \\ 
& \quad + (-1)^k \bigr( q^2 \hat{\sigma}_{(\overline{i}=\overline{1})
  (\overline{\mu} = \overline{k})} \hat{\sigma}_{(\nu=k-1)
  (\overline{i} = \overline{1})} + q^{-1} \hat{\sigma}_{(\nu=k-1)
  (\overline{i} = \overline{1})} \hat{\sigma}_{(\overline{i}=\overline{1})
  (\overline{\mu} =\overline{k})} \bigr) \\ 
&= (-1)^k (q+q^{-1}) \hat{\sigma}_{(\nu=k-1)(\overline{i} = \overline{1})}
  \hat{\sigma}_{(\overline{i}=\overline{1})(\overline{\mu} =\overline{k})} 
  + q \hat{\sigma}_{(\overline{i}=\overline{1})(\overline{\mu} = \overline{k})}
  \hat{\sigma}_{(\nu=k-1)(\overline{i} =\overline{1})} \bigr)\\ &= RHS
\end{align*}

\noindent as required.  In a similar way it can be shown that
equations \eqref{q2} and \eqref{q3} arise from the defining relations
of the $\hat{\sigma}_{ba}$.  Hence these compact defining relations
for the $\hat{\sigma}_{ba}$ incorporate not only the standard
$q$-Serre relations for the raising generators, but also the extra
ones, even though the equivalent $q$-Serre relations are not used in
the derivation.


\section{The $R$-matrix for the Vector Representation} \label{vector}

\noindent The Lax operator can be used to explicitly calculate an
$R$-matrix for any tensor product representation $\pi \otimes \phi$,
where $\phi$ is an arbitrary representation.  In particular, it
provides a more straightforward method of calculating $R$ for the
tensor product of the vector representation, $\pi \otimes \pi$, than
explicitly calculating projection operators \cite{Mehta}.

 By specifically constructing the $R$-matrix for the vector
 representation, we also illustrate concretely the way the recursion
 relations can be applied to find the $R$-matrix for an arbitrary
 representation.  Although the values for $\hat{\sigma}_{ba}$ obtained
 will change for each representation, they can always be constructed
 by applying the same equations in the same order.  We could choose to
 use only the relations listed in the tables in the appendix, but
 using the general form of the inductive relations shortens and
 simplifies the process.  Only the method of calculation is included
 here.  The full calculations follow those given in \cite{me}.

 First the vector representation is applied to the simple operators
 given in Table \ref{fundval}, with the results written below in Table
 \ref{basevector}.

\begin{table}[ht] 
\caption{The simple values of $\hat{\sigma}_{ba}$ in the vector
representation.}\label{basevector}
\centering
\begin{tabular}{|l|lll|} \hline
\multicolumn{1}{|c|}{Simple Root}& \multicolumn{3}{c|}{Corresponding
  $\hat{\sigma}_{ba}$} \\ \hline $\alpha_\mu = \delta_\mu -
  \delta_{\mu+1},\,\mu < k$ & $\hat{\sigma}_ {\mu(\mu+1)} $&$ =
  \hat{\sigma}_{(\overline{\mu+1})\overline{\mu}} $&$= E^\mu_{\mu+1}
  + E^{\overline{\mu+1}}_{\overline{\mu}}$ \\ $\alpha_s = \delta_k -
  \varepsilon_1,$ & $\hat{\sigma}_{(\mu=k)(i=1)} $&$ = (-1)^k q
  \hat{\sigma}_{(\overline{i}=\overline{1})(\overline{\mu} =
  \overline{k})} $&$= E^{\mu=k}_{i=1} + (-1)^k q E^{\overline{i} =
  \overline{1}}_{\overline{\mu} = \overline{k}}$ \\ $\alpha_t =
  \delta_k + \varepsilon_1, $&$ \hat{\sigma}_{(\mu=k)
  (\overline{i}=\overline{1})} $ &$ = (-1)^k q \hat{\sigma}_{(i=1)
  (\overline{\mu}=\overline{k})} $&$ =
  E^{\mu=k}_{\overline{i}=\overline{1}} + (-1)^k q
  E^{i=1}_{\overline{\mu}=\overline{k}} $\\ \hline
\end{tabular}
\end{table}

\noindent Then the inductive relations \eqref{indrel} are applied to
these simple operators to find the remaining values of
$\hat{\sigma}_{ba}$.  One of many equivalent ways of doing this is
given below.

 Construct
\begin{enumerate}
\item $\hat{\sigma}_{\nu\mu},\hat{\sigma}_{\overline{\mu}\,
  \overline{\nu}}$ for $1 \leq \nu < \mu \leq k$, using
  $\hat{\sigma}_{\mu (\mu+1)},\hat{\sigma}_{(\overline{\mu+1})
  \overline{\mu}}$ for $1 \leq \mu < k$ \label{item numu}
\item $\hat{\sigma}_{\mu (i=1)},
  \hat{\sigma}_{(\overline{i}=\overline{1}) \overline{\mu}}$ for $1
  \leq \mu \leq k$, using $\hat{\sigma}_{(\mu=k)(i=1)}$ and
  $\hat{\sigma}_{(\overline{i}=\overline{1}) (\overline{\mu}
  =\overline{k})} $ together with the values calculated in step
  \ref{item numu} \label{item mui}
\item $\hat{\sigma}_{(i=1)\overline{\mu}}, \hat{\sigma}_{\mu\,
  (\overline{i}=\overline{1})}$ for $1 \leq \mu \leq k$, using
  $\hat{\sigma}_{(\mu=k)(\overline{i}=\overline{1})}$ and
  $\hat{\sigma}_{(i=1)(\overline{\mu} =\overline{k})}$ and the results
  from step \ref{item numu}\label{item imubar}
\item $\hat{\sigma}_{\mu\,\overline{\nu}}$ for $1 \leq \mu,\nu \leq
  k$, using the results from steps \ref{item mui} and \ref{item
  imubar}
\end{enumerate}

\noindent Following this procedure in the vector representation, we
find the $R$-matrix for the vector representation of
$U_q[osp(2|n)]$, $\mathfrak{R} = (\pi \otimes \pi) \mathcal{R}$, is
given by

\begin{equation*}
\mathfrak{R} = q ^{\underset{a}{\sum} \pi(\tilde{h}_{a}) \otimes 
  \pi(\tilde{h}^{a})} \Bigl[ I \otimes I + (q -q^{-1}) \sum_{\varepsilon_{b} > 
  \varepsilon_{a}} (-1)^{[b]} E^a_b \otimes \hat{\sigma}_{ba} \Bigr]
\end{equation*}

\noindent where

\begin{equation*}
\hat{\sigma}_{ba} = q^{-(\varepsilon_a, \varepsilon_b)} E^b_a -
(-1)^{[b]([a]+[b])} \xi_a \xi_b q^{(\varepsilon_a, \varepsilon_a)}
q^{(\rho, \varepsilon_a - \varepsilon_b)}
E^{\overline{a}}_{\overline{b}}.
\end{equation*}

\noindent Writing this in a different form, we have the following result:

\begin{proposition}  The $R$-matrix for the vector representation, 
 $\mathfrak{R} = (\pi \otimes \pi) \mathcal{R}$, is given by 

\begin{equation*}
\mathfrak{R} = \sum_{a,b} q^{(\varepsilon_{a}, \varepsilon_{b})} E^a_a
\otimes E^b_b \; + \; (q-q^{-1}) \sum_{\varepsilon_{b} >
\varepsilon_{a}} (-1)^{[b]} E^a_b \otimes \tilde{\sigma}_{ba},
\end{equation*} 

\noindent where

\begin{equation*} 
\tilde{\sigma}_{ba} = E^b_a - (-1)^{[b]([a]+[b])} \xi_a \xi_b
q^{(\rho, \varepsilon_a - \varepsilon_b)}
E^{\overline{a}}_{\overline{b}}, \qquad \varepsilon_b > \varepsilon_a.
\end{equation*}
\end{proposition}

 We can also explicitly find the opposite $R$-matrix $\mathfrak{R}^T$
 as given in equation \eqref{RT}, using

\begin{equation*}
(E^a_b)^\dagger = (-1)^{[a]([a]+[b])} E^b_a.
\end{equation*}

\noindent We obtain this result:

\begin{proposition}
\noindent The opposite $R$-matrix for the vector representation,
$\mathfrak{R}^T = (\pi \otimes \pi) \mathcal{R}^T$, is given by

\begin{equation*}
\mathfrak{R}^T = \sum_{a,b} q^{(\varepsilon_a, \varepsilon_b)} E^a_a
\otimes E^b_b + (q-q^{-1}) \sum_{\varepsilon_b > \varepsilon_a}
(-1)^{[a]} E^b_a \otimes \tilde{\sigma}_{ab},
\end{equation*}

\noindent where

\begin{equation*}
\tilde{\sigma}_{ab} = E^a_b - (-1)^{[a]([a]+[b])} \xi_a \xi_b
  q^{(\rho, \varepsilon_a - \varepsilon_b)}
  E^{\overline{b}}_{\overline{a}}.
\end{equation*}

\end{proposition}

Although the root system was different, these formulae for
$\mathfrak{R}$ and $\mathfrak{R}^T$ on the vector representation take
the same form as those for $U_q[osp(m|n)],\, m >2,$ derived in \cite{us}.

\noindent{\bf Acknowledgements --} We gratefully acknowledge financial
support from the Australian Research Council.

\appendix

\section{The Relations Governing the Operators $\hat{\sigma}_{ba}$}

\begin{table}[h]
\caption{The relations for the operators $\hat{\sigma}_{ba}$ common to
all values of $m$} \label{list} \centering
\begin{tabular}{|ll|} \hline
\hspace{25mm} Relation & \hspace{10mm} Conditions \\ \hline

$\hat{\sigma}_{\nu(\mu+1)} = \hat{\sigma}_{\nu \mu}
  \hat{\sigma}_{\mu (\mu+1)} - q \hat{\sigma}_{\mu(\mu+1)}
  \hat{\sigma}_{\nu \mu}, $&$ \nu < \mu<k$ \\
$\hat{\sigma}_{(\overline{\mu+1}) \overline{\nu}} =
  \hat{\sigma}_{(\overline {\mu+1}) \overline{\mu}}
  \hat{\sigma}_{\overline{\mu} \overline{\nu}} - q
  \hat{\sigma}_{\overline {\mu} \overline{\nu}}
  \hat{\sigma}_{(\overline{\mu+1}) \overline{\mu}}, $&$ \nu < \mu < k$\\ 
$\hat{\sigma}_{b \overline{\mu}} = q^{(\alpha_\mu, \varepsilon_b)}
  \hat{\sigma} _{b(\overline{\mu+1})} \hat{\sigma}_{(\overline{\mu+1})
  \overline{\mu}} - q \hat{\sigma}_{(\overline{\mu+1}) \overline{\mu}}
  \hat{\sigma}_{b (\overline {\mu+1})}, $&$ \mu<k,\: \varepsilon_b >
  -\delta_{\mu + 1},$\\&$b\neq (\mu +1)$\\ 
$\hat{\sigma}_{\mu a} = q^{-(\alpha_\mu, \varepsilon_a)} 
  \hat{\sigma}_{\mu(\mu+1)} \hat{\sigma}_{(\mu+1)a} - q \hat{\sigma}_{(\mu+1)a}
  \hat{\sigma}_{\mu (\mu+1)},$&$ \mu<k,\: \varepsilon_a < \delta_{\mu +1},$\\ 
  &$a \neq (\overline{\mu + 1})$ \\ 
$\hat{\sigma}_{(\mu+1) \overline{\mu}} - \hat{\sigma}_{\mu(\overline{\mu+1})} 
  = q [\hat{\sigma}_{(\mu +1) (\overline{\mu+1})}, \hat{\sigma}_{\mu(\mu+1})],$
  & $\mu<k$ \\
$q^{(\alpha_\mu,\varepsilon_b)} \hat{\sigma}_{ba} \hat{\sigma}_{\mu (\mu +1)} 
  - q^{-(\alpha_\mu,\varepsilon_a)} \hat{\sigma}_{\mu (\mu+1)} 
  \hat{\sigma}_ {ba} = 0, $&$\mu<k;\: \varepsilon_b >\varepsilon_a;$ \\ 
  &$a\neq \mu, (\overline{\mu+1})$ \\ &and $b \neq (\mu+1), \overline{\mu}$ \\
$\hat{\sigma}_{\nu(i=1)} = \hat{\sigma}_{\nu (\mu=k)}
  \hat{\sigma}_{(\mu=k)(i=1)} - q \hat{\sigma}_{(\mu=k)(i=1)}
  \hat{\sigma}_{\nu (\mu=k)},$&$\nu <k$\\
$\hat{\sigma}_{(i=\overline{1})\overline{\nu}} =
  \hat{\sigma}_{(i=\overline{1}) (\overline{\mu} = \overline{k})}
  \hat{\sigma}_{(\overline{\mu} = \overline{k}) \overline{\nu}} - q
  \hat{\sigma}_{(\mu = \overline{k})\overline{\nu}} \hat{\sigma}_{(i =
  \overline{1})(\mu = \overline{k})}, $&$ \nu < k $\\
$\hat{\sigma}_{(\mu=k) \overline{\nu}} = q^{(\delta_\nu,\delta_k}
  \hat{\sigma} _{(\mu=k)(i=1)} \hat{\sigma}_{(i=1)\overline{\nu}} +
  \hat{\sigma}_{(i=1) \overline{\nu}} \hat{\sigma}_{(\mu=k)(i=1)},$&$
  \nu \leq k $\\ 
$\hat{\sigma}_{\nu (\overline{\mu} = \overline{k})} =q^{(\delta_\nu,\delta_k)} 
  \hat{\sigma}_{\nu (i= \overline{1})} 
  \hat{\sigma}_{(i = \overline{1})(\overline{\mu} = \overline{k})} + q^{-1} 
  \hat{\sigma}_{(i=\overline{1}) (\overline{\mu} =\overline{k})} 
  \hat{\sigma}_{\nu(i= \overline{1})}, $&$ \nu \leq k$\\ 
$\hat{\sigma}_{(\mu=k)(\overline{i} = \overline{1})} - (-1)^k q
  \hat{\sigma}_ {(i=1)(\overline{\mu} = \overline{k})} = q^{-1} [
  \hat{\sigma}_{(\mu=k)(i=1)},
  \hat{\sigma}_{(i=1)(\overline{i}=\overline{1})}],$& \\ 
$q^{(\alpha_s, \varepsilon_b)} \hat{\sigma}_{ba} \hat{\sigma}_{(\mu=k)(i=1)} 
  - (-1)^{[a]+[b]} q^{-(\alpha_s,\varepsilon_a)} \hat{\sigma}_{(\mu=k)(i=1)}
  \hat{\sigma}_{ba} = 0 $&$ \varepsilon_b > \varepsilon_a;\:
  \varepsilon_a \neq \delta_k, -\varepsilon_1$ \\ &and $\varepsilon_b
  \neq \varepsilon_1, -\delta_k$ \\
$\hat{\sigma}_{\nu(\overline{i}=\overline{1})} =
  \hat{\sigma}_{\nu(\mu=k)} \hat{\sigma}_{(\mu=k)(\overline{i}=\overline{1})} 
  - q \hat{\sigma}_{(\mu=k)(\overline{i}=\overline{1})} 
  \hat{\sigma}_{\nu(\mu=k)}, $&$ \nu<k$ \\
$\hat{\sigma}_{(\mu=k)\overline{\nu}} = q^{(\delta_k, \delta_\nu)}
  \hat{\sigma}_{(\mu=k)(\overline{i}=\overline{1})}
  \hat{\sigma}_{(\overline{i}=\overline{1})\overline{\nu}} + q^{-1}
  \hat{\sigma}_{(\overline{i}=\overline{1})\overline{\nu}}
  \hat{\sigma}_{(\mu=k)(\overline{i}=\overline{1})}, $&$\nu \leq k$ \\
$\hat{\sigma}_{(i=1) \overline{\nu}} = \hat{\sigma}_{(i=1)
  (\overline{\mu}= \overline{k})} \hat{\sigma}_{(\overline{\mu}=\overline{k}) 
  \overline{\nu}} - q\hat{\sigma}_{(\overline{\mu}=\overline{k})\overline{\nu}}
  \hat{\sigma}_{(i=1)(\overline{\mu}=\overline{k})}, $&$ \nu < k$\\
$\hat{\sigma}_{\nu (\overline{\mu}=\overline{k})} = q^{(\delta_k,\delta_\nu)} 
  \hat{\sigma}_{\nu(i=1)} \hat{\sigma}_{(i=1)(\overline{\mu}=\overline{k})} 
  + q^{-1} \hat{\sigma}_{(i=1)(\overline{\mu}=\overline{k})} 
  \hat{\sigma}_{\nu (i=1)}, $&$ \nu \leq k$ \\ 
$q^{(\alpha_t, \varepsilon_b)} \hat{\sigma}_{ba} \hat{\sigma}_{(\mu=k)
  (\overline{i}=\overline{1})} - (-1)^{[a]+[b]} q^{-(\alpha_t, \varepsilon_a)} 
  \hat{\sigma}_{(\mu=k)(\overline{i}=\overline{1})} \hat{\sigma}_{ba} = 0,$&$ 
  \varepsilon_b > \varepsilon_a;\;\varepsilon_a \neq \varepsilon_1, \delta_k$\\
  &$ \text{and }\varepsilon_b \neq -\varepsilon_1, -\delta_k$ \\ \hline
\end{tabular}
\end{table}


\end{document}